\documentclass{article}
\usepackage{graphicx}      
\usepackage{natbib}        

\usepackage{array}
\usepackage{amsmath,amssymb}

\newtheorem{Conj}{\bf{Conjecture}}

\def\PP{{\mathcal{P}}}
\def\UU{{\mathcal{U}}}
\def\TT{{\mathcal{T}}}

\title{Column generation for the discrete UC 
problem with min-stop ramping constraints} 

\author{Nicolas Dupin} 

\date{LRI, Universit\'e Paris-Sud, Universit\'e Paris-Saclay, France}

\begin{document}
\maketitle



\begin{abstract}                
The discrete unit commitment problem with  min-stop ramping constraints optimizes the daily production of thermal power plants. 
For this problem, compact Integer Linear Programming (ILP) formulations have been designed to solve exactly small instances and heuristically real-size instances.
This paper investigates whether Dantzig-Wolfe reformulation allows to improve the previous exact method and matheuristics.
The extended ILP formulation is presented with the column generation algorithm to solve its linear relaxation.
The experimental results show that the
Dantzig-Wolfe reformulation does not improve the quality of the linear relaxation of the tightest compact ILP formulations.
Computational experiments suggest also a conjecture which would explain such result: the compact ILP formulation of min-stop ramping constraints would be tight.
Such results validate the quality of the exact methods and matheuristics based on compact ILP formulations previously designed.
\end{abstract}



\section{Introduction}


Energy management induces complex production problems as electricity is not storable on a large scale.
This means that a large volume of electricity needs to be generated exactly at the time of consumption. 
However, power stations are not always able to keep up with the fluctuating demand. 
Dealing with power production and demand
induces several levels of optimization problems from strategic decisions in an uncertain environment  to daily production decisions (\cite{Ren93}).
\emph{Unit Commitment} (UC) problems denote these optimization problems, 
providing electricity according to the power demands and the power generation constraints
while minimizing the cost of the power generation. 
 
This paper  focuses on  a short term UC problem within a time window of two days and using a 30 minute discrete time step.
Power and reserves are generated with a  fleet of  coal, gas and fuel units.
The modulation capacity of such  thermal fleet is limited 
with ramping constraints (\cite{frangioni2008solving}).
In this paper, we consider a previously examined model by \cite{dupin2016tighter},  
the discretized UC problem with minimum stop and ramping constraints (UCPd) for thermal units.

 \cite{dupin2016tighter} provided several compact Integer Linear Programming (ILP) formulations for the UCPd problem,
 tightening the ILP formulations of min-stop ramping constraints
to have the best possible resolution using straightforwardly ILP solvers.
The resulting dual bounds of the Linear Programming (LP) relaxation are of good quality, 
providing gaps to the best known solutions in order of $1\%$ for the real size instances.
Matheuristics of \cite{dupin2018parallel} are based on the previous formulations and  solve efficiently large size instances within the short time limits imposed by the operational process.
This paper examines whether an extended Dantzig-Wolfe (D-W) reformulation can improve significantly the dual bounds.
Expected improvements allows two applications.
On one hand, a first issue is the acceleration of the exact methods.
On the other hands, using the extended formulation as a basis for matheuristics is another perspective.


%

This paper is organized as follows.
In section 2, we describe precisely the constraints of the  UCPd problem.
In section 3, we discuss related state-of-the-art elements. 
In section 4, we present a compact ILP formulation.
In section 5, the D-W reformulation of the previous 
ILP formulation is investigated with its column generation scheme to generate the linear relaxation.
The computational results are  presented in section 6, discussing theoretical and practical  implications of these results.

\begin{table}[ht]
 
       \centering
       \caption{Notation for UCPd}\label{notations}
\begin{tabular}{ll}
$u \in \mathcal{U}$ & index and set to designate  generating units. \\
 $t \in \mathcal{T}$ & index and set for  optimization time steps\\
    $i \in \mathcal{L}_u$  & index and set for operating points of unit $u$\\
  $N_u$  & Number of operating points for unit $u$.\\ 
  $P_{u,i}$  & Power generated by unit $u$ at point $i$. \\
  $R^1_{u,i}$   & Capacity in primary reserve for unit $u$ at point $i$.\\
  $R^2_{u,i}$ & Capacity in secondary reserve for unit $u$ at point $i$.\\
  $\Delta^{off}_{u}$ & Minimum down time for unit $u$. \\
  $\Delta^{on}_{u}$  & Minimum up time for unit $u$. \\
  $\Delta_{u,i}^+$ &  Min stop time at $i$ for unit $u$  before ramping up to $i+1$. \\
  $\Delta_{u,i}^-$ &  Min stop at  $i$ for unit $u$  before ramping down to $i-1$. \\
   $D_t^P$ & (Forecast) demand in power for period $t$. \\
   $D_t^{R1}$ & Demand in primary reserve for period $t$. \\
   $D_t^{R2}$ & Demand in secondary reserve for period $t$. \\
   $C_u^S$  & Start-up cost for unit  $u$. \\
   $C_u^F$ & Set-up  cost whenever unit $u$ is online.\\ 
   $C_u^P$  & Proportional cost to the power generated by unit $u$.\\
\end{tabular}
\end{table}

\section{Problem description}

This section presents the constraints of  UCPd.
 We refer to Table \ref{notations} for the notation.

\subsection{{Thermal UC with set-up and start-up costs}}
Basic UC decisions indicate the set-up status of the generators for each time period.
Production decisions are then assigned  to online generators  fulfilling power demands $D_t^{P}$  at any time step $t\in \mathcal{T}$.

A simple thermal UC problem can be formulated in Mixed Integer Linear Programming (MILP), considering that  units $u \in \mathcal{U}$  generate independently  power in  continuous domains $[P_u^{min}, P_u^{max}]$, minimizing the summed 
 operational costs for all the units:
 start-up costs $C_u^{S}$, set-up costs $C_u^{F}$, and  proportional costs $C_u^{P}$ to the power productions.
The decision variables are the power generated, $p_{u,t}\geqslant 0$, 
and the binary variables $x_{u,t},y_{u,t}\in \{0,1\}$ denoting respectively the \textit{set-up variables} and the \textit{start-up variables}. 
We have $x_{u,t}=1$ if and only if the unit $u$ is online at period $t$, whereas  $y_{u,t}=1$ indicates that the unit $u$ starts up at period $t$. 
In the following MILP, we consider furthermore 
\textit{min-up/min-down} constraints, which impose for all units $u \in \mathcal{U}$  minimal durations online $\Delta_u^{on}$  and  offline   $\Delta_u^{off}$:


\begin{eqnarray}
\displaystyle \min_{x,y\in\{0,1\}^M,P\geqslant 0} &\displaystyle\sum_{u \in \mathcal{U}}\sum_{t\in T} \left(C_u^{P}  p_{u,t} 
+  C_u^{F}    x_{u,t} +  C_u^{S} y_{u,t} \right) &  \label{objUCP}\\
 \forall u \in \mathcal{U},  \forall t \in \mathcal{T},  & x_{u,t} - x_{u,t-1}  \leqslant y_{u,t} &\label{couplUCP}\\
 \forall u \in \mathcal{U},  \forall t \in \mathcal{T}, &  P_u^{min} x_{u,t} \leqslant  p_{u,t}   &   \label{bornesProdUCP}\\
 \forall u \in \mathcal{U},  \forall t \in \mathcal{T}, &    p_{u,t}  \leqslant P_{u}^{max}  x_{u,t} &   \label{bornesProdUCP2}\\
\forall t \in \mathcal{T}, &  \sum_{u} p_{u,t}   = D_t^{P}  &   \label{demandUCP}\\
\forall u \in \mathcal{U}, \forall t \in \mathcal{T}, &  \sum_{t'=t-\Delta^{on}_{u}+1}^{t} y_{u,t'}  \leqslant x_{u,t} &  \label{minUpTurnOn} \\
\forall u \in \mathcal{U} ,\forall t \in \mathcal{T},& \sum_{t'=t-\Delta^{off}_{u}+1}^{t} y_{u,t'}  \leqslant 1- x_{u,t-\Delta^{off}_{u}+1} &  \label{minDownTurnOn}
\end{eqnarray}
The objective function (\ref{objUCP})  gathers start-up costs, set-up costs and  proportional costs  to the generated power,
it is linear once  variables $x,y,P$ are defined. 
Equation (\ref{couplUCP}) links the start-up variables to the set-up variables.
Equations (\ref{bornesProdUCP}) and (\ref{bornesProdUCP2}) bound the production domains when units are online, i.e. $x_{u,t}=1$, and impose zero production when $x_{u,t}=0$.
The productions match exactly the demands  at any time step with equation (\ref{demandUCP}). 
Equations (\ref{minUpTurnOn}) and (\ref{minDownTurnOn}) are the formulation of \textit{min-up/min-down} constraints from \cite{Tak05}.


\subsection{{Specific constraints for UCPd}}
The production domain is discrete for UCPd, power is generated only on \emph{operating points}
$i \in \mathcal{I}_u$ defined for each unit $u$.
The power associated to the operational point $i$ is $P_{u,i}$.  
The demands in power are not related to the discretization. 
To face such difficulties, the demand constraints are inequalities: over-productions are allowed. 
The minimization of the production costs dissuades  to over-generate. 
Mathematically, the production demands are equivalent to knapsack constraints  for each time period $t \in T$.
We consider also two  types or reserve differing in the operating delays (namely primary and secondary reserves).
Reserve constraints are modeled similarly to the power demands:  
defining for each operating point a maximal reserve participation, 
 the planning must  fulfill  reserve demands  at any time. 

\begin{figure}[ht]
      \centering
      \includegraphics[angle=0, width=.99\linewidth]{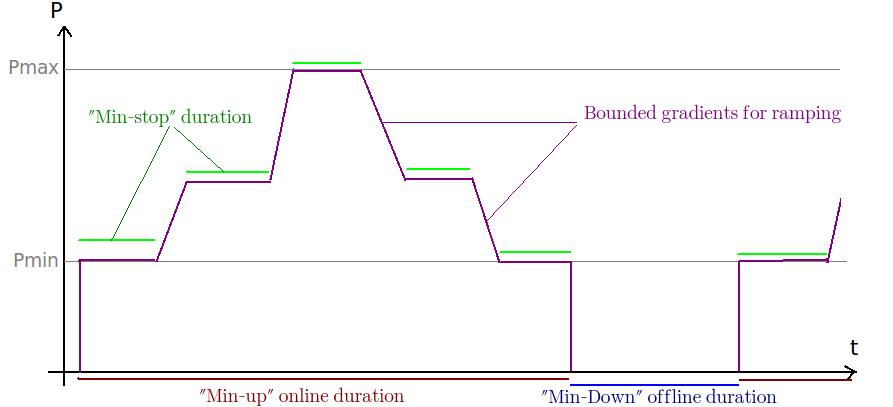}
	\caption{Illustration of the constraints for UCPd}\label{defConstrUCPd}
\end{figure}

\begin{figure}[ht]
      \centering
      \includegraphics[angle=0, width=.899\linewidth]{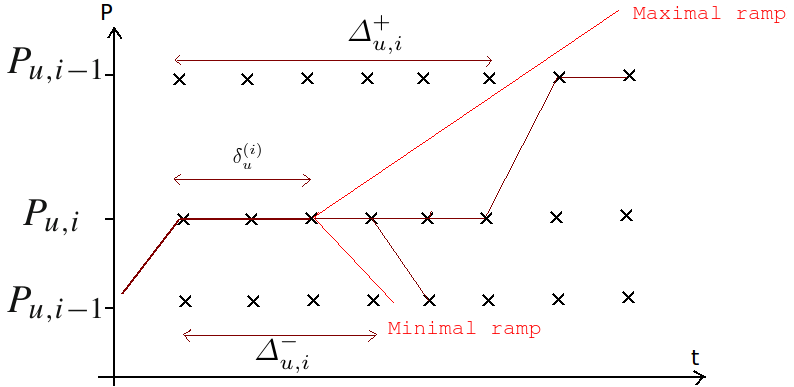}
	\caption{Illustration of the min-stop ramping constraints
	on the discretized operating points}\label{defDiscretUCPd}
\end{figure}

 
 To model possible power variations, three types of dynamic constraints are considered  and illustrated Figure \ref{defConstrUCPd}:

\begin{itemize}
 \item \textit{Min-up/min-down} constraints:  every unit $u$ has a minimum up time $\Delta^{on}_{u}$ online and a minimum down 
time $\Delta^{off}_{u}$ offline similarly to (\ref{minUpTurnOn}) and (\ref{minDownTurnOn}).
 \item \textit{Transition} constraints: 
When a unit generates at an operating point $i$ at period $t$, the allowable transitions for period $t+1$ are either to keep operating at point $i$
 or to shift to a neighboring point $j\in\{i-1,i+1\}$. 
 \item \textit{Min-stop ramping} constraints on operating points: 
Once unit $u$ produces on point $i$, the power must be stabilized during $\Delta_{u,i}^+$  (resp. $\Delta_{u,i}^-$) time steps
before reaching  point $i+1$  (resp. $i-1$), 
as illustrated in Figure~\ref{defDiscretUCPd}. 
\end{itemize}

\section{Related work}

ILP formulation results exist for \textit{min-up/min-down} constraints.
\cite{Tak00} provided a weak linear  formulation  with only set-up variables $x_{u,t}$,
\cite{Lee03} provided an exponential number of cuts to describe the convex hull of feasible integer points $x_{u,t}$ 
with a separation algorithm  for a Branch\&Cut implementation.
\cite{Tak05} proved that (\ref{minUpTurnOn}) and (\ref{minDownTurnOn}) dominate the previous formulation and cuts 
 and that the polytope defined with (\ref{couplUCP}),(\ref{minUpTurnOn}),(\ref{minDownTurnOn}) and $0 \leqslant x_{u,t},y_{u,t} \leqslant 1$ has integer extreme points.
\cite{Tak05} showed in experimental results  that the straightforward Branch\&Bound resolution
with (\ref{minUpTurnOn}) and (\ref{minDownTurnOn}) outperforms the Branch\&Cut algorithm derived from  \cite{Lee03,Tak00}.

A few decades ago, the solving capabilities of MILP solvers did not allow to consider 
 realistic models of UC problems. 
 MILP was widely used to model simple UC problems. 
With recent progress in the performances of  computers and MILP solving, more realistic UC models are considered.
Many works deal with additional  dynamic constraints on thermal production.
\cite{Arr06,morales2015tight} and \cite{gentile2017tight} provided efficient
MILP formulation to model start-up and shut-down trajectories when a thermal  production is in the domain  $[0,P_u^{min}]$.
\cite{silbernagl2016improving} and \cite{brandenberg2017tight}
 provided  efficient models to compute start-up costs and curves.
 \cite{frangioni2008solving} presented  different types and formulations for \emph{ramping} constraints 
when a thermal  production is in the domain $[P_u^{min},P_u^{max}]$, ensuring physical modulation constraints.
\cite{frangioni2009tighter} presented MILP formulations for the ramping constraints. 
These formulations have been strengthened with the addition of start-up and shut-down variables by \cite{ostrowski2012tight} and \cite{damci2013polyhedral}.
\cite{correa2017dynamic} provided recently more realistic formulations of ramping constraints.

UC problems consider mostly a continuous production domain. 
The discretization  is considered for the French case study, \cite{Dub05}
solved a short term UC   considering the whole French fleet with  a Lagrangian approach dualizing  demand constraints.
Each thermal unit induces a  sub-problem after  dualization.
The discretization  allows to solve these sub-problems independently with dynamic programming.
\cite{kruber2018resource} added constraints in the dynamic programming algorithm to solve thermal sub-problems.
 We  note that the influence of the discretization of ramping constraints
and the relaxation transition phases was studied by \cite{morales2017hidden}.

 Formerly, MILP models of UC problems were commonly solved  with Lagrangian decomposition dualizing demand constraints.
\cite{Liu00} and \cite{Dub05} developed to Lagrangian heuristics. 
D-W decomposition dualizing demand constraints are similar and allow to compute Lagrangian bounds.
\cite{fu2005long} and \cite{Roz12} investigated  D-W decomposition and Column Generation (CG) for long term UC problems.


\section{Compact ILP formulation}

Several variants to define  the variables of UCPd can be considered. 
Several compact formulations of some constraints are also possible once variables are defined.  
\cite{dupin2016tighter}  used projections and isomorphisms
to transform and compare  the polyhedrons defined by compact ILP formulations.
We present below one of the tightest ILP formulation provided by this work.


The production decisions are modeled using \emph{state variables} $s_{u,t}^{(i)} \in \{0,1\}$ are defined with 
$s_{u,t}^{(i)}=1$ if and only if the unit $u \in\mathcal{U}$ operates exactly at the point $i \in\mathcal{I}_u$ and period $t \in\mathcal{T}$, else $s_{u,t}^{(i)}=0$.

To have the tightest compact and linear formulation, additional variables are considered similarly to \cite{Tak05}.
\emph{Start-up variables} $y_{u,t}^{(i)-},y_{u,t}^{(i)+} \in \{0,1\}$
are  defined  for all $(u,t)\in\mathcal{U} \times \mathcal{T}$, and $i \in \mathcal{I}_u \cup \{0\}$ to indicate
if unit $u$ is ramped up  (resp. ramped down) to point $i$ at time $t$ from point $i-1$ (resp. from point $i+1$) at time $t-1$. 
\textit{Start-up variables}  are related with \textit{state variables}:  $ y_{u,t}^{(i)+}=s_{u,t}^{(i)} . s_{u,t-1}^{(i-1)}$, $ y_{u,t}^{(i)-}=s_{u,t}^{(i)} . s_{u,t-1}^{(i+1)}$.

To simplify the presentation of constraints, we extend 
the notations with  $y_{u,t}^{(0)+}=y_{u,t}^{(N_u-)}=0$ for all $t$.
The initial conditions are also considered with variables $s_{u,t}^{(i)},y_{u,t}^{(i)-},y_{u,t}^{(i)+}$ for $t\leqslant 0$ coding the previous production levels and moves.

It leads to the following ILP formulation: 

\begin{eqnarray}
 \displaystyle\min_{s,y} & \displaystyle\sum_{i} \left(C_u^{P} P(u,i) + C_u^{F}\right)  s_{u,t}^{(i)} 
+ \sum_{u,t} C_u^{S} y_{u,t}^{(1)+} \\
\forall u,t,    &  \sum_i  s_{u,t}^{(i)}  \leqslant 1  \label{gub}\\ 
\forall u,t,i,  & y_{u,t}^{(i-1)-} + \sum_{j\geqslant i} \left(s_{u,t}^{(j)} - s_{u,t-1}^{(j)} \right)  = y_{u,t}^{(i)+} \label{couplEtat}\\
\forall u,t,i,  &  \sum_{j\geqslant i} s_{u,t}^{(j)} \leqslant  \sum_{j\geqslant i-1} s_{u,t+1}^{(j)}\label{transEtat1}\\
\forall u,t,i,  & \sum_{j\geqslant i} s_{u,t}^{(j)} \geqslant \sum_{j\geqslant i+1} s_{u,t+1}^{(j)} \label{transEtat2}\\
\forall u,t,i,  & \sum_{t'=t+1}^{t+\Delta_{u,i}^+} y_{u,t'}^{(i+1)+}  + 
 \sum_{t'=t+1}^{t+\Delta_{u,i}^-} y_{u,t'}^{(i-1)-}   \leqslant s_{u,t}^{(i)}\label{minStop}\\
\forall u,t , & \sum_{t'=t-\Delta_{u}^{on}+1}^{t} y_{u,t'}^{(1)+}  \leqslant \sum_{i>0} s_{u,t}^{(i)}\label{Tak1}\\
\forall u,t, & \sum_{t'=t-\Delta_{u}^{off}+1}^{t} y_{u,t'}^{(1)+}  \leqslant 1- \sum_{i>0} s_{u,t-\Delta_{u}^{off}+1}^{(i)}\label{Tak2}\\
\forall t,  &   \sum_{u,i} P(u,i) \phantom{1}  s_{u,t}^{(i)}   \geqslant  D_t^{P} \label{DemandState}\\
 \forall t,  & \sum_{u,i} R_1(u,i) \phantom{1}  s_{u,t}^{(i)}  \geqslant  D_t^{R1} \label{R1State}\\
\forall t,  & \sum_{u, i} R_2(u,i) \phantom{1}  s_{u,t}^{(i)}   \geqslant  D_t^{R2} \label{R2State}
\end{eqnarray}

Constraints (\ref{gub}) are implied by the definition of variables $s$: for all time period $t$, a unit $u$ produces at one operating point $i>0$ or is offline.
Constraints (\ref{couplEtat}) linearize the coupling constraints $ y_{u,t}^{(i)+}=s_{u,t}^{(i)} . s_{u,t-1}^{(i-1)}$, $ y_{u,t}^{(i)-}=s_{u,t}^{(i)} . s_{u,t-1}^{(i+1)}$.
The tightest formulation of {transition} constraints and  the \textit{min-stop ramping} constraints. are respectively (\ref{transEtat1})-(\ref{transEtat2}) and (\ref{minStop}), we refer to \cite{dupin2016tighter}.
Constraints (\ref{Tak1}) and (\ref{Tak2}) are the \textit{min-up/min-down} constraints
similarly to (\ref{minUpTurnOn}) and (\ref{minDownTurnOn}), noticing that
the set-up and start-up variables $x_{u,t},y_{u,t}$ are  $x_{u,t} =  \sum_{i>0} s_{u,t}^{(i)}$ and $y_{u,t}=y_{u,t}^{(1)+}$. 
Constraints (\ref{DemandState}-\ref{R2State}) are demand constraints in generated power, primary and secondary reserves. 



\section{Extended ILP formulation} 

This section investigates the D-W reformulation of previous ILP model, dualizing the demand constraints in power and reserves (\cite{Van00}).
The extended ILP formulation is presented with the \emph{Column Generation} (CG) algorithm  to solve its LP relaxation.


\subsection{Extended ILP formulation}

We denote by $\PP^u$  the set of all the feasible production planning for unit $u$
considering the initial conditions and the technical constraints (\ref{gub}-\ref{Tak2}).

The variables of the extended formulation are indexed  by $\PP= \bigcup_u\PP^u$, 
the set of all the possible production plannings.
Binaries $z_p \in  \{0,1\}$ are equal to $1$ if and only if the production planning $p \in \PP$ is used in the global production planning.
For each $p \in \PP$, we denote by
 $P_{p,t}$, $R^1_{p,t}$, $R^2_{p,t}$ the production and reserves generated in
 the planning $p$ at period $t$. 
 It gives rise to the following ILP formulation:

\begin{eqnarray}
\min_{z_p \in \{0,1\}} &\sum_{p \in \PP} c_p z_p& \label{objExt}\\
\forall t \in \TT, \phantom{2}   &  \sum_{p \in \PP} P_{p,t}z_p \geqslant D_t^{P}  \label{demUCPdExt}\\
\forall t \in \TT, \phantom{2}   &  \sum_{p \in \PP} R^1_{p,t}z_p \geqslant D_t^{R1} \\
\forall t \in \TT, \phantom{2}   &  \sum_{p \in \PP} R^2_{p,t}z_p \geqslant D_t^{R2}  \label{R2UCPdExt}\\
\forall u \in \UU,  \phantom{2}  & \sum_{p \in \PP_{u}} z_p = 1 & \label{convexUCPd}
\end{eqnarray}

Constraints (\ref{demUCPdExt}-\ref{R2UCPdExt}) are the demands in generated power and reserve capacities.
Constraints (\ref{convexUCPd}) is required to express that a single production planning
is assigned to each unit.

The dynamic constraints are induced in the definition of production patters $\PP$.
It is polyhedrally equivalent to  consider the convex hull of the integer points defined by
constraints (\ref{gub}-\ref{Tak2}). 
The LP relaxation of the extended formulation
furnishes Lagrangian bounds, these dual bounds are at least as good as the  LP relaxation of the compact ILP formulations (we refer to \cite{Van00}).

The difficulty to deal with this formulation is that the set of feasible
plannings has an exponential size and can be enumerated only for very small instances.
Hence, CG techniques are required to deal with a reasonable number of variables.


\subsection{Column Generation scheme}

The LP relaxation of (\ref{objExt}-\ref{convexUCPd}) is solved by the CG algorithm, adding iteratively new production patterns.
Having a subset of production patters $\PP^o \subset \PP$,  the  \emph{Restricted Master Problem} (RMP) denotes the LP relaxation of the previous ILP formulation restricted to the variables
indexed by $\PP^o $.
Defining $\PP_{u}^o = \PP_u \cap \PP^o$, the RMP is written as following with the dual variables:

\begin{equation}
  \begin{array}{lll}
\mbox{RMP}(\PP^o)& =\displaystyle \min_{z_p \geqslant 0} \displaystyle\sum_{p \in \PP^o} c_p z_p \\
\forall t \in \TT,    &  \displaystyle\sum_{p \in \PP^o} P_{p,t}z_p \geqslant D_t^{P} & \left(\pi^P\right) \\
\forall t \in \TT,    &  \displaystyle\sum_{p \in \PP^o} R^1_{p,t}z_p \geqslant D_t^{R1} & \left(\pi^{R1}\right) \\
\forall t \in \TT,    &  \displaystyle\sum_{p \in \PP^o} R^2_{p,t}z_p \geqslant D_t^{R2} & \left(\pi^{R2}\right) \\
\forall u \in \UU,    & \displaystyle\sum_{p \in \PP_{u}^o} z_p = 1 & (\sigma) \\
\forall p \in \PP^o,  & z_p \geqslant 0
  \end{array}
\end{equation}

The constraints (\ref{convexUCPd}) imply $z_p \leqslant 1$
and thus the LP relaxation can be written using only positivity constraints $z_p \geqslant 0$.
As is, there are no dual variables associated to constraints $z_p \leqslant 1$.

The reduced cost related to the variable $z_p$, denoted by $\mbox{RC}(z_p)$, has following value :
\begin{equation*}
 \mbox{RC}(z_p) = c_p + \sigma_u - \sum_t \left(\pi_t^P P_{p,t} + \pi_t^{R1} R^1_{p,t} + \pi_t^{R2} R^2_{p,t}\right)
\end{equation*}
where $u$ is such that $p \in \PP_u$.


The CG algorithm iterates while there exist no columns $p \in \PP$ such that $ \mbox{RC}(z_p)<0$.
In this case, $\mbox{RMP}(\PP^o)$ is the value of the LP relaxation of (\ref{objExt}-\ref{convexUCPd}).
Otherwise, variables with a negative reduced cost are added in the RMP
and the procedure is repeated till the stopping criterion is reached.

Solving $\mbox{RC}^* = \min_{p \in \PP^o} \mbox{RC}(z_p)$, the CG sub-problems, contains
the two previous computations: $\mbox{RC}^*\geqslant 0$ is equivalent to the termination criterion,
and optimal solutions have a negative reduced cost when $\mbox{RC}^*< 0$ and can be selected to add in the RMP for the next iteration.

Computations to optimality of  $\mbox{RC}^*$ are useful only for the last iteration, to prove the termination of the CG algorithm.
Otherwise, heuristics are useful to generate quickly columns with a negative reduced cost
and also to spend less time solving sub-problems.

We note that the CG algorithm requires the feasibility of $\mbox{RMP}(\PP^o)$ at each iteration for the
computations of the  dual variables for the following sub-problems. 
The CG algorithm shall be initialized with columns ensuring the feasibility. 
The feasibility is ensured for each iteration adding columns.
 Removing columns with a null value in the last RMP  has no effect in the feasibility of the next RMP, it allows to deal with bounded sizes
 of RMP.
A first initialization strategy is to consider for each unit its maximal and minimal production planning regarding the technical constraints and the initial conditions.
Otherwise, a global feasible planning can be computed quickly with the heuristics described by \cite{dupin2016matheuristics}.

\subsection{Solving CG subproblems}

The key point is to solve efficiently $\mbox{RC}^* = \min_{p \in \PP} \mbox{RC}(z_p)$.
The optimization are decomposed for each unit, $\mbox{RC}^* = \min_{u} \mbox{RC}^*_u$,
where $\mbox{RC}^*_u= \min_{p \in \PP_u} \mbox{RC}(z_p)$ is the minimization problem over the feasible production patterns for unit $u$.
Computations of $\mbox{RC}_u^*$ are independent and can be computed in parallel.
It can be formulated using the compact ILP formulations of Section 4:

\begin{eqnarray}
  \mbox{RC}_u^*=& \min_{x,y}  f_{cout}^u(x,y)+ f_{dual}^u(x)\\ 
\forall t,    &  \sum_i  x_{t}^{(i)}  \leqslant 1 \label{SPeq1} \\ 
\forall t,i,  &  \sum_{j\geqslant i} x_{t}^{(j)} \leqslant  \sum_{j\geqslant i-1} x_{t+1}^{(j)}\\
\forall t,i,  & \sum_{j\geqslant i} x_{t}^{(j)} \geqslant \sum_{j\geqslant i+1} x_{t+1}^{(j)} \\
\forall t,i,  & y_{t}^{(i-1)-} + \sum_{j\geqslant i} \left(x_{t}^{(j)} - x_{t-1}^{(j)} \right)  = y_{t}^{(i)+}\label{couplXY}\\
\forall t,i,  & \sum_{t'=t+1}^{t+\Delta_{u,i}^+} y_{t'}^{(i+1)+}  + 
 \sum_{t'=t+1}^{t+\Delta_{u,i}^-} y_{t'}^{(i-1)-}   \leqslant x_{t}^{(i)}\\
\forall t , & \sum_{t'=t-\Delta_{u}^{on}+1}^{t} y_{t'}^{(1)+}  \leqslant \sum_{i>0} x_{t}^{(i)}\\
\forall t, & \sum_{t'=t-\Delta_{u}^{off}+1}^{t} y_{t'}^{(1)+}  \leqslant 1- \sum_{i>0} x_{t-\Delta_{u}^{off}+1}^{(i)}\label{SPeq2}\\
\forall t,i,  &y_{t}^{(i)+},y_{t}^{(i)-},x_{t}^{(i} \in \{0,1\}
\end{eqnarray}

where the objective function is composed of $f_{cost}^u(x,y)$, the production costs like in the compact formulation, and $f_{dual}^u(x)$ are additional dual costs:

$\displaystyle f_{dual}^u(x) = \sigma_u  - \sum_{u,i,t} \left( \pi^P_t  P_{u,i} + \pi^{R1}_t  R^1_{u,i} + \pi^{R2}_t  R^2_{u,i} \right) \phantom{1}  x_{t}^{(i)} $


$\displaystyle f_{cost}^u(x,y) = \sum_{t} \left( C^{F}_u  x_{t}^{(1)}
+  C^{S}_u y_{t}^{(1)+}
+ \sum_{i} C^{P}_u  P_{u,i}  x_{t}^{(i)} \right)
$

If $\mbox{RC}_u^*$ is written as an ILP, it can be computed with a dynamic programming algorithm,
ensuring a polynomial complexity to the sub-problem resolution.

\section{Experimental results}

The computational experiments used the instance dataset from \cite{dupin2016tighter}.
These instances were generated from real-world data for the French thermal fleet.
Time horizon is $2$ days, with $96$ periods of $30$ minutes. 
There are up to $80$ generating units, with around $3$ discrete points per unit.

 Tests were computed with an   Intel(R) Core(TM) i5-4430 CPU, 3.00GHz, running Linux, with 4 CPU cores.
ILP and LP were  solved with Cplex 12.6 using the OPL interface as a first implementation. %
A first issue is to compare the quality of the LP relaxation of the extended formulation and the compact formulations from \cite{dupin2016tighter},
and to examine if LP relaxation improvements have a positive impact on the exact resolution.
A second issue is to derive primal heuristics from  the CG algorithm.

It is a known result that the LP relaxation with D-W reformulation is at least as good as the ones of compact formulations.
More precisely, the D-W reformulation is equivalent to consider the
convex hull of the polyhedron defined by the subproblems.
The equality of both LP relaxations are reached when the polyhedrons defined by the sub-problems have integer extreme points
in the compact formulation (\cite{Van00}).
The improvements of the LP relaxation provided by  \cite{dupin2016tighter}  also
tightened the sub-problems defined for each unit, which was a first way to close the gap between the first LP relaxation
to the LP relaxation with a convexification of sub-problems.
Actually, the LP relaxation of sections 4 and 5 have exactly the same value for each considered instance.
These empirical results on the LP relaxation suggest to conjecture that the tightest compact formulations
provide tight formulations for the min-stop ramping constraints:

\vskip 0.4cm

\begin{Conj}
The polytope defined for a single unit with constraints (\ref{SPeq1})-(\ref{SPeq2}) 
 describes the convex hull of the feasible integer points for the min-stop ramping constraints. 
 In other words, the ILP formulation of constraints (\ref{SPeq1})-(\ref{SPeq2}) is tight, and the Dantzig-Wolfe reformulation of UCPd 
dualizing demand and reserve constraints has the same LP relaxation than the compact formulation of section 4 and the other
tightest and equivalent formulation derived by  \cite{dupin2016tighter}.
\end{Conj}
 \vskip 0.2cm

Such conjecture generalizes the polyhedral results obtained by \cite{Tak05}.
This conjecture was verified computing sub-problems (\ref{SPeq1})-(\ref{SPeq2}) as LP or ILP 
during CG iterations. 
Finding a counter-example where the LP relaxation and  the ILP resolution of a CG sub-problem have different optimal values, 
this would  induce that the conjecture is false. This was never observed, opening the perspective to prove formally the conjecture.

\begin{figure}[ht]
      \centering
      \includegraphics[angle=0, width=.78996\linewidth]{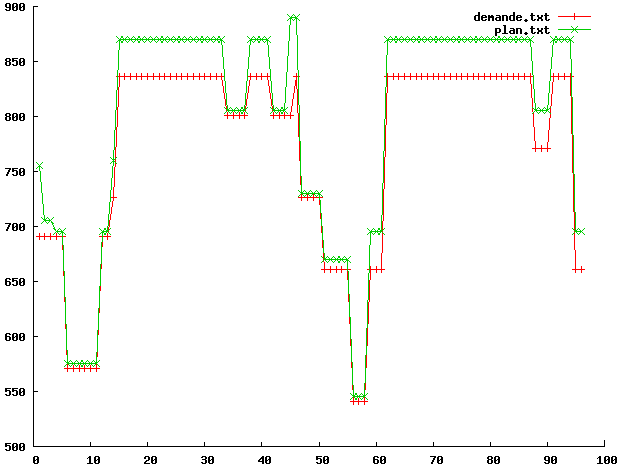}
	\caption{Demand covering with constraints (\ref{DemandState})-(\ref{R2State})}\label{solutionUCPd}
\end{figure}

We note that the convergence of the CG scheme is difficult, requiring much more time than compact LP relaxations.
The dual variables are very unstable, with a erratic convergence.
Indeed, for most of the iterations, the dual variables indexed by the time periods have few non-zeros.
It tends to generate columns with peaks of production in the periods with non-zeros, and null productions where it is possible with the min-stop ramping and min-up/min-down constraints.
This structure of dual variables is actually implied by the inequalities (\ref{DemandState}-\ref{R2State}).
Power discretization and the min-stop ramping constraints imply that few constraints (\ref{DemandState}-\ref{R2State})
reach the lower bounds imposed by the constraints as illustrated in Figure \ref{solutionUCPd}.
The complementary slackness theorem ensures that the corresponding dual variables of the non-saturated constraints are null. 
Such types of production planning given by the CG iterations are not efficient to be combined in global and feasible integer solutions.
This closes the perspectives to use CG iterations for a primal heuristic based on the integer RMP and CG iterations, and
validates also the matheuristics from \cite{dupin2018parallel}.

These results validate to use in practice compact ILP formulations for UCPd,
in a Branch\&Bound resolution and in the computation of LP relaxation for the LP-based heuristics from \cite{dupin2018parallel}.
These results for UCPd are  specific to this simple UC model. 
Modeling constraints (\ref{DemandState}) as inequalities induced difficulties for the convergence of the CG algorithm,
whereas the minimization of the slackness to the power demands like in Goal Programming approaches would be more convenient for a CG scheme and CG heuristics.
Lastly, a gap between compact LP relaxations and convexified extended formulations
shall be analyzed considering in the UC model other types of power plants and/or more realistic constraints for thermal power plants, like in the work of \cite{rottner2018combinatorial}.

\end{document}